\documentclass[12pt]{article}
\usepackage{amsfonts}
\usepackage{enumerate}
\usepackage{amsmath,amssymb}
\usepackage{latexsym}

\usepackage{float}
\restylefloat{table}

\newtheorem{thm}{Theorem}

\newtheorem{note}{Note}

\newtheorem{lema}{Lemma}

\makeatother

\begin{document}

\title
{\bf Classification of 
quasi-symmetric 2-$(64,24,46)$ designs of Blokhuis-Haemers type }
\author
{Dean Crnkovi\'c\\
{\it\small Department of Mathematics} \\
{\it\small University of Rijeka} \\
{\it\small Radmile Matej\v ci\'c 2, 51000 Rijeka, Croatia}\\[5pt]
B. G. Rodrigues \\
{\it\small School of Mathematics, Statistics and Computer Science} \\
{\it\small University of KwaZulu-Natal} \\
{\it\small Durban 4041, South Africa}\\[5pt]
Sanja Rukavina \\
{\it\small Department of Mathematics} \\
{\it\small University of Rijeka} \\
{\it\small Radmile Matej\v ci\'c 2, 51000 Rijeka, Croatia}\\[5pt]
and \\[5pt]
Vladimir D. Tonchev\\
{\it\small Department of Mathematical Sciences,} \\
{\it\small Michigan Technological University, Houghton, MI 49931, USA}\\[-15pt] 
}
\date{}
\maketitle
\begin{abstract}

This paper completes the classification of quasi-symmetric 2-$(64,24,46)$ designs
 of Blokhuis-Haemers type
supported by the dual code $C^{\perp}$ of the binary linear code $C$ spanned by the
lines of $AG(3,2^2)$ initiated in \cite{bgr-vdt}. It is shown that $C^{\perp}$ contains
exactly 30,264 nonisomorphic quasi-symmetric 2-$(64,24,46)$ designs 
obtainable from maximal arcs in $AG(2,2^2)$ via the Blokhuis-Haemers construction. 
The related strongly regular graphs are also discussed.

\end{abstract}

\section{Introduction} 

We assume  familiarity with basic facts, terminology and notation
 from  design theory, finite geometry, and coding theory
  \cite{ak, bjl, CRC,  Hir,  vdt}. 
For strongly regular graphs, cf. \cite{aeb-g}, \cite{BCN}, and
for the theory of quasi-symmetric designs, one may consult the monograph \cite{MS}.

\medskip

In \cite{b-h}, Blokhuis and Haemers  
gave an elegant construction of a quasi-symmetric
design $D(q)$ with parameters
2-$(q^3,q^2(q-1)/2,q(q^3-q^2-2)/4)$ and block intersection
numbers $q^2(q-2)/4$ and $q^2(q-1)/4$,  where $q$ is an arbitrary power of 2.
The Blokhuis-Haemers construction is a clever refinement of a
method due to Shrikhande and Raghavarao \cite{s-r}.
Every block of $D(q)$ is the union
of $q(q-1)/2$ parallel lines in the 3-dimensional affine geometry $AG(3,q)$
labeled by a  block of a symmetric
2-$(q^2, q(q-1)/2, q(q-2)/4)$ design invariant under the translation group of $AG(2,q)$
and defined in terms of maximal arcs in $AG(2,q)$.

It  a recent paper \cite{max-arcs},  Jungnickel and Tonchev studied the properties
of $D(q)$ and proved that the
number of nonisomorphic quasi-symmetric designs obtainable 
via the Blokhuis-Haemers construction grows exponentially with
linear growth of $q$.
Following \cite{max-arcs},
we call  any design obtained via the Blokhuis-Haemers construction a BH-design.

In particular, if $q=4$, it was proved in \cite[Proposition 3.7]{max-arcs} that there are at least
28,844 nonisomorphic 2-$(64,24,46)$ BH-designs.
In this paper, we give a complete classification of
 all  2-$(64,24,46)$ BH-designs up to isomorphism, and show that the total
number of nonisomorphic such designs is exactly 30,264 (Theorem \ref{th4}).
This classification is a continuation of the work in
\cite{bgr-vdt}. Our approach employs a binary linear code 
 associated with the  designs in question
that utilizes the following property of BH-designs:
if $q>2$, every block of $D(q)$ meets every line of $AG(3,q)$ in
an even number of points \cite{b-h}, \cite[Lemma 3.2, (d)]{max-arcs}.
This  property implies the following.

\begin{lema}
\label{l1}
\cite{bgr-vdt}.
If $q > 2$,  every block of ${D}(q)$ is the support of a codeword of
weight $q^2(q-1)/2$ in the dual code $C^{\perp}$ of the binary code $C$ of length $q^3$ spanned
by the incidence vectors of the lines in $AG(3,q)$.
\end{lema}
  
If $q=4$,  the binary code $C$ spanned by the lines of $AG(3,4)$
is of dimension 51 (by Hamada's rank formula \cite{hamada}),
hence the dimension of $C^{\perp}$ is 13.  
The weight enumerator $W(x)$ of $C^{\perp}$ is  
\[
\label{w}
W(x)=1 + 1008x^{24} + 6174x^{32} + 1008x^{40} + x^{64}.
\]  
The automorphism group $G=Aut(C^{\perp})$ of  $C^{\perp}$ coincides with the collineation group $\Gamma L(3,4)$
of $AG(3,4)$, and is of order 
\begin{equation}
\label{G}
 23,224,320 = 2 \times 4^3(4^3 -1)(4^3-4)(4^3-4^2) = 2^{13}\cdot 3^4 \cdot 5 \cdot 7. 
\end{equation}

We use the code $C^{\perp}$ to find block-by-point incidence matrices
of quasi-symmetric 2-$(64,24,46)$ BH-designs as collections of 336 codewords
of weight 24, such that every two codewords share either 8 or 12 nonzero positions.
For these computations, we used Magma \cite{magma}  and Cliquer \cite{cliquer}.

In \cite{bgr-vdt}, this approach was used to classify up to isomorphism 
 all BH-designs invariant under
automorphisms of $C^{\perp}$ of odd prime order. It was shown that there is exactly 
one isomorphism class of designs admitting automorphisms of order 7, fifteen isomorphism classes of
designs admitting automorphisms of order 5, and no designs with automorphisms of order 3. In addition,
it was shown that there is exactly one BH-design in $C^{\perp}$ 
with a full automorphism group being the Sylow 2-subgroup of $G$  of order $2^{13}$. 

It is the goal of this paper to complete the classification of 2-$(64,24,46)$ BH-designs,  by
finding representatives of the isomorphism
classes of all remaining designs, which must have full automorphism groups of
 order $2^i$ for $i< 13$.

\section{Counting 2-$(64,24,46)$ BH-designs}
\label{sec2}

The following lemmas were crucial in making the classification
of 2-$(64,24,46)$ BH-designs computationally feasible.

\begin{lema} 
\label{21}
The total number of distinct 2-$(64,24,16)$ BH-designs is  $3^{21}$.
\end{lema}

{\bf Proof}.  We consider two designs $D_1$, $D_2$ to be distinct if their collections of blocks
are distinct, that is, there is a block $B_1$ of $D_1$ which is not a block of $D_2$. 
The statement of the lemma follows from \cite[Theorem 2.2]{max-arcs}, and can also be verified
computationally as follows.

Let $D'$ be the 2-$(64,4,1)$ design of the lines in $AG(3,4)$, and consider the
natural resolution of $D'$ into 21 parallel classes of lines, where each
parallel class consists of a line through the origin and its translates.
Let $P$ be such a parallel class, consisting of lines $L_1, \ldots, L_{16}$,
and let $D''$ be a symmetric 2-$(16,6,2)$ design such that the points of $D''$
are labeled with the sixteen lines from $P$. The Blokhuis-Haemers construction
(or more generally, the Shrikhande-Raghavarao construction \cite{s-r})
replaces the 16 lines of $P$ with 16 new blocks of size 24,
each being a union 
\[ U=L_{i_1}\cup L_{i_2} \cup \cdots \cup L_{i_6} \]
of six lines of $P$ that correspond to a block of $D''$.
Clearly, any two of the sixteen new blocks of size 24
share two lines from $P$, hence meet each other in 8 points. If $D''$ is appropriately
chosen design obtained from maximal arcs in $AG(2,4)$ (cf. \cite{max-arcs}), 
the resulting new design
with blocks of size 24 is a quasi-symmetric 2-$(64,24,46)$ design  \cite{b-h}, \cite{max-arcs}.

By Lemma \ref{l1}, every block of a BH-design meets every line
of $AG(3,4)$ in an even number (0, 2 or 4) of points. It is easy to check by computer
that there are exactly 48 unions, $U_1, U_2, \cdots, U_{48}$,
of six lines from $P$ that meet every line of $AG(3,4)$ evenly.
We define a graph $\Gamma_P$ with vertices $U_1, U_2, \cdots, U_{48}$, where
$U_i$ and $U_j$ are adjacent if they share either 8 or 12 points.
A quick check shows that the maximum clique size in $\Gamma_P$ is 16,
and there are exactly three 16-cliques. The three 16-cliques
 are mutually disjoint, that is,
partition the vertex set $\{ U_i \}_{i=1}^{48}$. In addition, every 16-clique
consists of blocks that meet each other in exactly 8 points. We call the three 16-cliques
associated with a parallel class P {\it special} cliques.
The symmetric 2-$(16,6,2)$ design $D''$ associated with any special
clique is invariant under the elementary abelian group $E_{16}$
of order 16 acting transitively on the blocks of $D''$ and the set of lines of $P$. The blocks
of $D''$ are maximal arcs in an affine plane
of order 4, $AG(2,4)$, associated with $P$ \cite{max-arcs}.
The design $D''$ is isomorphic to the unique 2-$(16,6,2)$ design
admitting a 2-transitive automorphism
group,
 and is also the unique SDP design with these parameters, in the terminology of \cite{JT, K}.

Since
 the collineation group of $AG(3,4)$ is transitive on the set of 21 parallel classes of lines,
the same applies to each parallel class. It follows that the collection of blocks of any BH-design
is a union of 21 special cliques, one clique for each of the 21 parallel classes.
This implies that the number of BH-designs is at most $3^{21}$.

To prove the equality, we compute  a graph $\Delta$ having 1008 vertices corresponding
to the blocks of size 24 associated with the 21 parallel classes 
(48 blocks per parallel class). We define
two blocks to be adjacent if they share 8 or 12 points. A quick computer check shows that
any two blocks associated with different parallel classes are adjacent in $\Delta$.
Consequently, every collection of 21 special cliques, one for each of the 21 different
parallel classes, is the set of blocks of a BH-design. {$\Box$} 

\begin{note}
{\rm
The 1008 blocks of size 24 in Lemma \ref{21}
correspond to the 1008 codewords of weight 24 in the code $C^{\perp}$ (cf. (\ref{w})).
}
\end{note}

\begin{lema}
\label{l3}
Every 2-$(64,24,46)$ BH-design is invariant under an elementary abelian group 
of order 64, isomorphic to the translation group $T$ of $AG(3,4)$.
\end{lema}

{\bf Proof}.
 The translation group  $T$ stabilizes each parallel class $P$ of $AG(3,4)$ and the special cliques
associated with $P$. 
For a geometric proof of an analogous result for arbitrary $q$ see \cite{max-arcs}. $\Box$

Since the automorphism group
$G \cong \Gamma L(3,4)$ of $C^{\perp}$ preserves the partition of the set of 1008 codewords
into triples of special 16-cliques associated with the 21 parallel classes of lines in $AG(3,4)$,
we have the following.

\begin{lema}
\label{l4}
The automorphism group of every BH-design is a subgroup of
$G = Aut(C^{\perp}) \cong \Gamma L(3,4)$.
\end{lema}

Lemma \ref{l4} also follows from \cite[Lemma 3.4]{max-arcs}.

Suppose that the number of nonisomorphic BH-designs in $C^{\perp}$ is $N$, and let $D_1, 
\ldots, D_N$
be a set of $N$ pairwise nonisomorphic BH-designs. 
As a corollary of Lemmas \ref{21}, \ref{l3}, and \ref{l4}, we have the following equation:
\begin{equation}
\label{eq}
3^{21}= \sum_{i=1}^N \frac{|G|}{|Aut(D_i)|}, 
\end{equation}
where $Aut(D_i)$ denotes the full automorphism group of $D_i$.

We can split the right-hand side sum in (\ref{eq}) into two parts:
\begin{equation}
\label{eq2}
 \sum_{i=1}^N \frac{|G|}{|Aut(D_i)|} = \frac{|G|}{64}N_{64} +\sum_{j: \ 128 \big|  |Aut(D_j)|} \frac{|G|}{|Aut(D_j)|},
\end{equation} 
where $N_{64}$ denotes the number of nonisomorphic designs with full group of order 64.

Thus, if we find the number of nonisomorphic BH-designs having automorphism group of
order divisible by $2^7 = 128$, we can determine $N_{64}$ and  $N$
from equations (\ref{eq}) and (\ref{eq2}), which would complete the classification of
2-$(64,24,16)$ designs.

\section{Classifying BH-designs with a group of order 128}

It is  known that every finite group of order $2^n$ contains a
subgroup of order $2^i$ for every $i$ in the range $1 \le i \le n$
(cf. \cite[Theorem 6.5, page 116]{rose}).
By this property, finding all nonisomorphic BH-designs in $C^{\perp}$
which
are invariant under a subgroup of $G=Aut(C^{\perp})$ of order 128,
will complete the
classification of BH-designs in  $C^{\perp}$.

By Lemma \ref{l3},  it suffices to
consider only the subgroups of $G$ which contain the translation group $T$ of order 64.
Using Magma, we found that the group $G$ contains 962 conjugacy classes of subgroups of 
order 128, but only two  contain representatives that contain $T$ as a subgroup.
In what follows, we will use two such subgroups of $G$, denoted by  $H_1$ and $H_2$.

The group $H_1$ is isomorphic to the group labeled by $(128,2163)$ in the Magma small groups library \cite{magma}.
The normalizer $N_G(H_1)$ of $H_1$ in $G$ has order 73728. We use the normalizer $N_G(H_1)$ 
for elimination of isomorphic designs.

The group $H_1$ partitions the 1008 codewords of $C^{\perp}$ of weight 24 into 39 orbits:
 15 orbits of length 16, and 24 orbits of length 32.
We call an orbit {\it good } if any two codewords from that orbit share exactly 8 or 12 nonzero positions. 
It turns out that all $H_1$-orbits are good, hence all orbits 
could be used to build a BH-design with parameters $2$-$(64,24,46)$.
We call two orbits {\it compatible} if every codeword
from one orbit shares
 exactly 8 or 12 nonzero positions with every codeword from the other orbit.
We define a graph $\Gamma$ with 39 vertices corresponding to the $H_1$-orbits,
two vertices being adjacent if and only if the corresponding orbits are compatible.
Every BH-design invariant under $H_1$ corresponds to a clique in $\Gamma$
labeled by a set of pairwise compatible orbits containing a total of
 336 codewords.
Let $\Gamma_1$ (resp. $\Gamma_2$)
 be the subgraph of $\Gamma$ having as vertices the 15 orbits of
length 16 (resp. the 24 orbits of length 32). The maximum clique size
in $\Gamma_1$ is 5, while the maximum clique size in $\Gamma_2$
is 8. It follows that any BH-design invariant under $H_1$ corresponds to the
 union of one 5-clique from $\Gamma_1$ and one 8-clique from $\Gamma_2$.
A further check shows that every $H_1$-orbit of length 16 is compatible 
with all $H_1$-orbits of length 32. 
Thus, the union of each 5-clique of $\Gamma_1$ with an 
8-clique  $\Gamma_2$ gives a BH-design.

The normalizer $N_G(H_1)$ acts in two orbits on the set of 15 $H_1$-orbits of length 16, 
while it acts transitively on the set of 24 $H_1$-orbits of length 32.
Hence, it is sufficient to consider only designs that contain one fixed $H_1$-orbit of length 32,
for example, the first such orbit. 
For further elimination of isomorphic designs, we use
the stabilizer of the first $H_1$-orbit of length 32 in the group $N_G(H_1)$.

A computation based on this approach shows that there are exactly 2688 mutually 
nonisomorphic BH-designs that admit $H_1$ as an automorphism group.
Information about the orders of the full automorphism groups of these designs,
and the number of nonisomorphic designs with a full automorphism group of given order,
is given in {\rm T}able \ref{table-designs-H1}.

\smallskip
\begin{center}
\begin{table}[htpb!]
\begin{footnotesize}
\begin{tabular}{r r | r r | r r}
 \hline
$|{\rm Aut}({\mathcal D})|$ & \# designs & $|{\rm Aut}({\mathcal D})|$ & \# designs & $|{\rm Aut}({\mathcal D})|$ & \# designs
 \cr \hline \hline
20480  & 2 & 1280 &  1 & 512 &    64 \\
8192   & 1 & 1024 &  8 & 256 &   210 \\
2048   & 3 &  640 & 12 & 128 &  2387 \\
\cr \hline \hline
\end{tabular}\caption{\footnotesize Nonisomorphic $2$-$(64,24,46)$ BH-designs admitting $H_1$ as an automorphism group} \label{table-designs-H1}
\end{footnotesize}
\end{table}
\end{center}

The second subgroup $H_2$ of order 128 which contains $T$, is isomorphic to the group labeled 
by $(128,1578)$ in the Magma small groups library. The normalizer $N_G(H_2)$ has order 21504.
$H_2$ partitions the 1008 codewords of weight 24 into 35 orbits: 
7 orbits of length 16, and 28 orbits of length 32.
All $H_2$-orbits of length 16 are good, and 21 of the 28 orbits of length 32 are good orbits.
As in the case with $H_1$, 
we define a graph $\Gamma$ with 28 vertices corresponding to good $H_2$-orbits,
two vertices being adjacent if and only if the corresponding orbits are compatible, 
and search for cliques that determine $2$-$(64,24,46)$ BH-designs.
 We define graphs $\Gamma_1$ and $\Gamma_2$ as induced subgraphs of $\Gamma$ with vertex sets 
determined by orbits of length 16 or 32, respectively.
$\Gamma_1$ is a complete graph on 7 vertices, and the maximum size of a clique in $\Gamma_2$ is 7. 
Any $H_2$-orbit of length 16 is compatible with all good $H_2$-orbits of length 32,
hence any 7-clique in $\Gamma_2$ 
together with the vertices of  $\Gamma_1$ determines a BH-design with parameters $2$-$(64,24,46)$. 
There are exactly 2187 7-cliques in
$\Gamma_2$,  yielding 2187 distinct $2$-$(64,24,46)$ designs.
 In that set of 2187 BH-designs, there are 17 mutually nonisomorphic ones.
Information about the orders of the full automorphism groups of these designs,
and the number of nonisomorphic designs with a full automorphism group of given order,
is listed in {\rm T}able \ref{table-designs-H2}.

\smallskip
\begin{center}
\begin{table}[htpb!]
\begin{footnotesize}
\begin{tabular}{r r | r r }
 \hline
$|{\rm Aut}({\mathcal D})|$ & \# designs & $|{\rm Aut}({\mathcal D})|$ & \# designs
 \cr \hline \hline
8192  & 1 & 512 &  3 \\
1024  & 1 & 256 &  1 \\
896   & 1 & 128 & 10 \\
\cr \hline \hline
\end{tabular}\caption{\footnotesize Nonisomorphic $2$-$(64,24,46)$ BH-designs admitting $H_2$ as an automorphism group} \label{table-designs-H2}
\end{footnotesize}
\end{table}
\end{center}

The full automorphism groups of order 8192, 1024, 512 and 256 contain
 both groups $H_1$ and $H_2$, hence among the 17 BH-designs invariant under $H_2$
there are 6 designs that admit $H_1$ as well.
Thus, there are exactly 2699 nonisomorphic designs among the $2$-$(64,24,46)$
BH-designs summarized  in {\rm T}able \ref{table-designs-H1}
 and {\rm T}able \ref{table-designs-H2}. 
We note that the $2$-$(64,24,46)$ BH-designs found in \cite{bgr-vdt}
have full automorphism groups of orders $20480=2^{12}\cdot 5$, $8192=2^{13}$,
 $1280=2^8\cdot 5$, $896=2^7\cdot 7$, and $640=2^7\cdot 5$.
Thus, all designs from \cite{bgr-vdt} admit also an automorphism group of order 128.
Hence, we have the following.

\begin{thm} \label{BH-designs-128}
There are exactly $2699$ nonisomorphic  $2$-$(64,24,46)$ BH-designs admitting 
an automorphism group of order $128$.
\end{thm}

%By the results of Section \ref{sec2}, 
If ${\mathcal D}$ is  a  $2$-$(64,24,46)$ BH-design that does not admit an automorphism group of order $128$,
then the translation group $T$ of $AG(3,4)$ is its full automorphism group.
 
Theorem \ref{BH-designs-128} and the results from Section \ref{sec2} imply the main result of this paper.  
\begin{thm}
\label{th4}
There are exactly $30,264$ nonisomorphic BH-designs with parameters $2$-$(64,24,46)$.
\end{thm}
\noindent {\bf Proof}.
Using the data from Tables \ref{table-designs-H1} and \ref{table-designs-H2},
we have
\begin{equation}
\label{eq4}
\sum_{j: \ 128 \big|  |Aut(D_j)|} \frac{|Aut(C^{\perp}|}{|Aut(D_j)|} =457566003.
\end{equation}
whence, from equations (\ref{eq}) and (\ref{eq2}) we have
\[ 3^{21}= \sum_{i=1}^N \frac{|Aut(C^{\perp}|}{|Aut({\mathcal D_i})|}=  \frac{23224320}{64}N_{64} + 457566003. \]
Hence, $N_{64}=27565$ and $N=30264$.
$\Box$

\bigskip

Information about the orders of the full automorphism
 groups of the 30264 nonisomorphic BH-designs,
and the number of nonisomorphic designs with a group of given order, 
is listed in {\rm T}able \ref{table-designs-all}.

\smallskip
\begin{center}
\begin{table}[htpb!]
\begin{footnotesize}
\begin{tabular}{r r | r r | r r}
 \hline
$|{\rm Aut}({\mathcal D})|$ & \# designs & $|{\rm Aut}({\mathcal D})|$ & \# designs & $|{\rm Aut}({\mathcal D})|$ & \# designs
 \cr \hline \hline
20480  & 2 & 1024 &  8 & 256 &   210 \\
8192   & 1 &  896 &  1 & 128 &  2397 \\
2048   & 3 &  640 & 12 &  64 & 27565 \\
1280   & 1 &  512 & 64 &     &       \\
\cr \hline \hline
\end{tabular}\caption{\footnotesize Nonisomorphic BH-designs with parameters $2$-$(64,24,46)$} \label{table-designs-all}
\end{footnotesize}
\end{table}
\end{center}

Among the 2699 nonisomorphic 
BH-designs with parameters 2-(64,24,46) admitting an automorphism group of order 128, there are
three designs having  2-rank equal to 12, namely the design with full automorphism group of order 8192, 
one of the designs with  full automorphism group of order 2048,  and one of the designs with  full automorphism 
group of order 20480. 
All other designs have  2-rank  13.
A list of the nonisomorphic  $2$-$(64,24,46)$ BH-designs admitting an automorphism group of order 128 
is available  at 

\begin{verbatim} 
www.math.uniri.hr/~sanjar/structures/.
\end{verbatim}  

\section{Strongly regular graphs with parameters (336,80,28,16)}

A {\it strongly regular graph} with parameters $(n,k,\lambda,\mu)$
is a graph with $n$ vertices which is regular of degree $k$,
every two adjacent vertices have $\lambda$
common neighbors, and every two nonadjacent
vertices have $\mu$ common neighbors.
The block graph of a quasi-symmetric 2-$(64,24,46)$ design
with block intersection numbers 8 and 12,
where two blocks are adjacent if they share 12 points, is
a strongly regular graph with parameters $(336,80,28,16)$.
These are also the parameters of the block graph of any
Steiner 2-$(64,4,1)$ design, where two blocks are adjacent if they share a point.
The graphs obtained from a quasi-symmetric  2-$(64,24,46)$ BH-design
and a resolvable Steiner 2-$(64,4,1)$ design
(for example, the 2-$(64,4,1)$ design of the lines of $AG(3,4)$)
share the property that their sets of vertices
can be partitioned into 21 cocliques of size 16.
Strongly regular graphs whose point set can be partitioned into cliques
(or cocliques), are studied in \cite{HT}.
A $(336,80,28,16)$ strongly regular graph which is not the block graph
of a Steiner or a quasi-symmetric design is discussed in \cite{aeb} and \cite{jenrich}
(see also
 the on-line table of strongly regular graphs maintained by Andries Brouwer \cite{aeb-g}).

The block graphs of the 2699 nonisomorphic  $2$-$(64,24,46)$  BH-designs 
admitting  an automorphism of order 128, split into 2371 
isomorphism classes of strongly regular graphs with parameters $(336,80,28,16)$. 

Information about orders of the full automorphism groups 
of these strongly regular graphs is given in {\rm T}able \ref{table-SRGs}.

\smallskip
\begin{center}
\begin{table}[htpb!]
\begin{footnotesize}
\begin{tabular}{r r | r r }
 \hline
$|{\rm Aut}(\Gamma)|$ & \# SRGs & $|{\rm Aut}(\Gamma)|$ & \# SRGs
 \cr \hline \hline
245760 & 1 & 1280 &    4 \\
61440  & 1 & 1024 &   12 \\
24576  & 1 &  896 &    1 \\
8192   & 1 &  640 &    7 \\
6144   & 1 &  512 &   56 \\
4096   & 2 &  384 &    4 \\
2048   & 6 &  256 &  220 \\
1536   & 3 &  128 & 2051 \\
\cr \hline \hline
\end{tabular}\caption{\footnotesize Nonisomorphic block graphs of $2$-$(64,24,46)$ BH-designs admitting an automorphism group of order 128} \label{table-SRGs}
\end{footnotesize}
\end{table}
\end{center}

\section{ Acknowledgements} 

The work of D. Crnkovi\'c and S. Rukavina has  been fully supported 
by the  Croatian Science Foundation under the Project 1637.
B. G. Rodrigues acknowledges research support by the National Research Foundation of South Africa 
(Grant Numbers 84470 and 91495). V. D. Tonchev  acknowledges support by NSA Grant H98230-16-1-0011.


\begin{thebibliography}{30}

 \bibitem{ak}
E. F. Assmus Jr., J. D. Key, {\it Designs and their Codes}, Cambridge University Press, Cambridge, 1992.

\bibitem{bjl}
T. Beth, D. Jungnickel, H. Lenz, {\it Design Theory}, 2nd Edition, Cambridge University Press, Cambridge, 1999.

\bibitem{b-h}
A. Blokhuis, W. H. Haemers, An infinite family of quasi-symmetric designs,
{\it J. Statist. Plann. Inference} {\bf 95} (2001), 117 -119.

\bibitem{magma}
W. Bosma, J. Cannon, {\it Handbook of Magma Functions},\\
Department of Mathematics, University of Sydney, November 1994,\\ http://magma.maths.usyd.edu.au/magma.


\bibitem{aeb-g} A. ~E. ~Brouwer, Parameters of Strongly Regular Graphs,\\ 
http://www.win.tue.nl/~aeb/graphs/srg/srgtab.html.

\bibitem{aeb}
A. ~E. ~Brouwer, A strongly regular graph on 336 vertices, preprint, September 12, 2014,
 http://www.win.tue.nl/~aeb/preprints/srg336.pdf.


\bibitem{BCN} A. E. Brouwer, A. M. Cohen, A. Neumaier, {\it Distance Regular Graphs}, Springer, 1989.

\bibitem{CRC} \bibitem{CRC} C. J. Colbourn, J. F. Dinitz, Eds., 
{\it Handbook of Combinatorial Designs},
Second Edition, Chapman \& Hall/CRC, 2007.
 

\bibitem{hamada}
N. Hamada, On the p-rank of the incidence matrix of a balanced or partially balanced incomplete
block design and its application to error correcting codes, {\it Hiroshima Math. J.} {\bf 3} (1973), 153-226.

\bibitem{HT} W. Haemers, V.D. Tonchev, Spreads in strongly regular graphs,
{\it Des. Codes Cryptogr.}  {\bf 8} (1996), 145-157. 

\bibitem{Hir} J.~W.~P. Hirschfeld: \emph{Projective Geometries over Finite Fields (2nd edition).} Oxford University Press (1998).
 

\bibitem{jenrich}
T. Jenrich, New strongly regular graphs derived from the $G_2(4)$ graph,
\texttt{arXiv:1409.3520v2}, 15 Sept. 2014.

\bibitem{JT} D. Jungnickel and V. D. Tonchev, On symmetric and quasi-symmetric
designs with the symmetric difference property and their codes,
{\it J. Combin. Theory}, Ser. A {\bf 59} (1992), 40-50.

\bibitem{max-arcs}
D. Jungnickel, V. D. Tonchev, Maximal arcs and quasi-symmetric designs,
{\it Des. Codes Cryptogr.} {\bf 77} (2015), 365--374.

\bibitem{K} W. M. Kantor, Symplectic groups, symmetric designs and line
ovals, {\it J. Algebra} {\bf 33} (1975), 43-58.

\bibitem{cliquer}
S. Niskanen, P. R. J. \"Osterg\r ard, Cliquer User's Guide, Version 1.0. Tech. Rep. T48,
Communications Laboratory, Helsinki University of Technology, Espoo, Finland, 2003.

%\bibitem{robinson}
%D. Robinson, {\it A Course in the Theory of Groups}, Springer-Verlag, New York Berlin Heidelberg, 1996.

\bibitem{bgr-vdt}
B. G. Rodrigues, V. D. Tonchev, On Quasi-symmetric 2-(64,24, 46)
Designs Derived from Codes, in: R. Pinto et al. (eds.), Coding Theory and Applications,
CIM Series in Mathematical Sciences 3, Springer International Publishing, Switzerland, 2015, pp. 327--333.

\bibitem{rose}
H. E. Rose, {\it A Course on Finite Groups}, Springer, 2009.

\bibitem{MS} M. S. Shrikhande, S. S. Sane, {\it Quasi-symmetric Designs}, Cambridge University Press,
Cambridge, 1991.

\bibitem{s-r}
S. S. Shrikhande, D. Raghavarao,
A method of construction of incomplete block designs, {\it Sankhy\={a}}, Ser. A {\bf 25} (1963), 399 -402.

\bibitem{vdt}
V. D. Tonchev,
{\it Combinatorial Configurations: Designs, Codes, Graphs},
Pitman Monographs and Surveys in Pure and Applied Mathematics 40,
Wiley, New York, 1988.


\end{thebibliography}
\end{document}